\documentclass[draft]{elsart}

\usepackage{amsmath, amssymb, amscd}

\def\frak{\mathfrak}
\def\Bbb{\mathbb}
\def\Cal{\mathcal}

\newcommand{\gr}{\operatorname{gr}}

\newcommand{\Ad}{\operatorname{Ad}}

\renewcommand{\ker}{\operatorname{ker}}
\newcommand{\im}{\operatorname{im}}

\newcommand{\Aut}{\operatorname{Aut}}

\newcommand{\Rho}{P}

\newcommand{\fg}{{\frak g}}
\newcommand{\tg}{{\tilde{\mathfrak g}}}
\newcommand{\tp}{{\tilde{\mathfrak p}}}
\newcommand{\x}{\times}

\let\ccdot\cdot
\def\cdot{\hbox to 2.5pt{\hss$\ccdot$\hss}}

\newcommand{\al}{\alpha}

\newcommand{\ga}{\gamma}

\newcommand{\om}{\omega}
\renewcommand{\phi}{\varphi}
\newcommand{\ph}{\varphi}
\newcommand{\ps}{\psi}
\newcommand{\si}{\sigma}

\newcommand{\ze}{\zeta}

\newcommand{\Ga}{\Gamma}
\newcommand{\La}{\Lambda}
\newcommand{\Om}{\Omega}
\newcommand{\Ph}{\Phi}
\newcommand{\Ps}{\Psi}

\begin{document}
\begin{frontmatter}
  \title{On Nurowski's Conformal Structure Associated to a Generic
    Rank Two Distribution in Dimension Five}

\author{Andreas \v Cap\thanksref{FWF}\corauthref{COR}}
\corauth[COR]{Corresponding author.}
\ead{Andreas.Cap@esi.ac.at}
\author{Katja Sagerschnig\thanksref{FWF}}
\address{Fakult\"at f\"ur Mathematik, Universit\"at Wien, Nordbergstra\ss
  e 15, A--1090 Wien, Austria}
\ead{Katja.Sagerschnig@univie.ac.at}
\thanks[FWF]{Both authors supported by project P 19500--N13 of the ``Fonds
  zur F\"orderung der wissenschaftlichen Forschung'' (FWF)}



\begin{abstract}
  For a generic distribution of rank two on a manifold $M$ of
  dimension five, we introduce the notion of a generalized contact
  form. To such a form we associate a generalized Reeb field and a
  partial connection. From these data, we explicitly constructed a
  pseudo--Riemannian metric on $M$ of split signature. We prove that a
  change of the generalized contact form only leads to a conformal
  rescaling of this metric, so the corresponding conformal class is
  intrinsic to the distribution.
  
  In the second part of the article, we relate this conformal class
  to the canonical Cartan connection associated to the distribution.
  This is used to prove that it coincides with the conformal class
  constructed by Nurowski.
\end{abstract}
\end{frontmatter}

\section{Introduction}\label{1}
The study of generic rank two subbundles in the tangent bundles of
five--dimensional manifolds goes back to Elie Cartan's famous ``five
variables paper'' \cite{Cartan:five} from 1910. This paper is
remarkable in several respects. First, by constructing a canonical
Cartan connection associated to such distributions, Cartan showed that
they have non--trivial local invariants. Second, for the simplest
instance of such a distribution, Cartan showed that the infinitesimal
symmetries form an exceptional simple Lie algebra of type $G_2$.  This
was the first instance of an exceptional simple Lie algebra showing up
``in real life''.

In modern terminology, the homogeneous model of generic rank two
distributions in dimension five is the quotient of the split real form
of an exceptional Lie group $G$ of type $G_2$ by a maximal parabolic
subgroup $P\subset G$. An explicit description of this homogeneous
space and its generic rank two distribution can be found in
\cite{Sa2}. Cartan's construction associates to an arbitrary generic
rank two distribution on a five--manifold $M$ a Cartan geometry of
type $(G,P)$.

In pioneering work culminating in \cite{Tanaka}, N.~Tanaka showed that
Cartan geometries modelled on the quotient of a semisimple group by a
parabolic subgroup are (under small restrictions, which have later
been eliminated) equivalent to simpler underlying geometric
structures.  These geometric structures have been intensively studied
during the last years under the name ``parabolic geometries''. The
description of the structures equivalent to parabolic geometries is
best phrased in terms of filtered manifolds, an up to date version can
be found in \cite{Srni05}. In particular, for many geometries this
underlying structure is only a filtration of the tangent bundle of a
manifold with certain properties, see \cite{Sa1}. In the special case
of the split real form of $G_2$ with the appropriate parabolic
subgroup, this recovers Cartan's result. 

For generic rank two distributions in dimension five, progress was
made recently by P.~Nurowski in \cite{Nurowski}.  Using the canonical
Cartan connection associated to a distribution $\Cal H$ on $M$,
Nurowski constructed a canonical conformal structure of signature
$(2,3)$ on $M$. Rather than a rank two distribution, Nurowski's
starting point was an underdetermined system of ODEs of certain type,
which can be equivalently described by a generic distribution. The
system of ODEs is determined by a single smooth function, and in his
article Nurowski gives an impressive (and frightening) formula for a
metric in the conformal class in terms of this function.

It was soon realized, see \cite{Srni05}, that Nurowski's construction
can be interpreted as an analog of the Fefferman construction, which
to a non--degenerate CR--manifold of hypersurface type associates a
canonical conformal structure on the total space of a circle bundle.
The main point in this analogy is not that a conformal structure
occurs in both cases, but the interpretation of the construction in
terms of Cartan connections. In terms of this analogy, Nurowski's
construction corresponds to the version of the Fefferman construction
based on Cartan connections, which was developed in \cite{BDS}.

The aim of this paper is to present a description of Nurowski's
conformal class which is analogous to J.~Lee's description of the
Fefferman construction, see \cite{Lee1,Lee2}. In section \ref{2}, we
introduce the notion of a generalized contact form for a generic rank
two distribution $\Cal H$ on a five--manifold $M$. Starting from such
a form $\al$ we explicitly construct a pseudo--Riemannian metric
$g_\al$ of signature $(2,3)$ on $M$. Then we prove that another choice
of generalized contact form leads to a conformally related metric, so
the conformal class $[g_\al]$ is intrinsic to the distribution $\Cal
H$. The formula for $g_\al$ in terms of data derived from $\al$ is
rather simple and explicit, and requires no knowledge of the canonical
Cartan connection.

In section \ref{3} we show that our conformal structure coincides with
the one constructed by Nurowski. This is based on the theory of
Weyl--structures for parabolic geometries, which interprets the
canonical Cartan connection in terms of underlying data. This also
explains how the formula for the metric $g_\al$ was actually found.
The developments in section \ref{3} are also interesting from the
point of view of the general theory of parabolic geometries, since
this is the first time that essential parts of a Weyl structure are
explicitly computed for one of the more involved parabolic
geometries. This article is based on results obtained during the work
on the second author's PhD thesis, which will contain a complete
description of this Weyl structure.

\section{A canonical conformal structure}\label{2}
In this section, we first introduce the notion of a generalized
contact form for a generic rank two distribution $\Cal H$ on a smooth
manifold $M$ of dimension five. Given a generalized contact form
$\al$, we construct a pseudo--Riemannian metric $g_\al$ on $M$. Then
we show that the metrics associated to different generalized contact
forms are always conformal to each other. Hence the conformal class
$[g_\al]$ on $M$ depends only on the distribution $\Cal H$.

\subsection{Generic rank $2$ distributions in dimension $5$}\label{2.1}
We start by collecting some facts on such distributions and fixing some
notation. Recall that a rank $2$ distribution $\mathcal{H}$ on a
$5$--manifold $M$ is called \textit{generic} if the values of linear
combinations of iterated Lie brackets of at most three sections of
$\Cal H$ in each $x\in M$ span the tangent space $T_xM$.  In
particular, sections of $\Cal H$ and their Lie brackets have to span a
rank $3$ subbundle $[\mathcal{H},\mathcal{H}]$. Defining
$T^{-1}M=\mathcal{H}$ and $T^{-2}M=[\mathcal{H},\mathcal{H}]$, we
obtain a filtration
$$
T^{-1}M\subset T^{-2}M\subset TM
$$
of the tangent bundle by smooth subbundles. We use the convention
that $T^iM=0$ for $i\geq 0$ and $T^iM=TM$ for $i\leq -3$. Then the
filtration is compatible with the Lie bracket of vector fields in the
sense that for $\xi\in\Gamma(T^iM)$ and $\eta\in\Gamma(T^jM)$ we get
$[\xi,\eta]\in \Gamma(T^{i+j}M)$. For $i=-1,-2,-3$ we define
$\gr_i(M)=T^iM/T^{i+1}M$ and then $\gr(M)=\oplus_{i=-3}^{-1}\gr_i(TM)$
is the associated graded vector bundle to the tangent bundle. For
$i=-2,-3$, we denote by $q_i:T^iM\to\gr_i(TM)$ the natural quotient
map.

The Lie bracket of vector fields induces a skew symmetric bilinear
bundle map $\{\ ,\ \}:\gr(TM)\times\gr(TM)\to\gr(TM)$, called the Levi
bracket, which is homogeneous of degree zero. Note that the two
nontrivial components
$$
\{\ ,\ \}:\Lambda^2 T^{-1}M\to
\gr_{-2}(TM);\quad\{\xi(x),\eta(x)\}:=q_{-2}([\xi,\eta](x)) 
$$
and 
$$
\{\ ,\ \}:T^{-1}M\otimes\gr_{-2}(TM)\to\gr_{-3}(TM);\quad
\{\xi(x),q_{-2}(\ze(x))\}:=q_{-3}([\xi,\ze](x))
$$
are both isomorphisms of vector bundles.

The filtration of $TM$ dualizes to filtration of the cotangent bundle
$T^*M$ and we consider the associated graded bundle
$$
\gr(T^*M)=\oplus_{i=1}^3\gr_i(T^*M).
$$
By construction, $\gr_i(T^*M)=(\gr_{-i}(TM))^*$.

\subsection{Generalized contact forms and Reeb fields}\label{2.2}
\begin{defn}\label{def22}
  Let $\Cal H\subset TM$ be a generic rank two distribution on a
  $5$--manifold $M$. A \textit{generalized contact form} for $\Cal H$
  is a smooth section $\al$ of the bundle $(T^{-2}M)^*$ such that for
  each $x\in M$ the kernel of the linear map $\al(x):T^{-2}_xM\to\Bbb
  R$ is $T^{-1}_xM=\Cal H_x$.
\end{defn}
By definition, a generalized contact form is a partially defined one
form. We will see below, that $\al$ can be canonically extended to a
true one form on $M$. Note the the condition on the kernel implies
that $\al$ is nowhere vanishing.

Given a generalized contact form $\al$ for $\Cal H$, we next want to
introduce an analog of the Reeb vector field. Consider a local section
$r$ of $T^{-2}M$ which is transversal to $T^{-1}M$, i.e.~such that
$\ph:=q_{-2}(r)\in\Ga(\gr_{-2}(TM))$ is nowhere vanishing. Then for
$\xi,\eta\in\Ga(T^{-1}M)$, there is a unique smooth section
$\nabla_\xi\eta\in\Ga(T^{-1}M)$ such that
\begin{equation}
  \label{nabla-1-def}
\{\nabla_\xi\eta,\ph\}=q_{-3}([\xi,[\eta,r]]).   
\end{equation}
By construction, the operator
$\nabla:\Ga(T^{-1}M)\times\Ga(T^{-1}M)\to\Ga(T^{-1}M)$ is bilinear
over $\Bbb R$. Using that $q_{-3}(\xi)=0$ and
$q_{-3}([\eta,r])=\{\eta,\ph\}$ one immediately concludes from the
defining equation that $\nabla_{f\xi}\eta=f\nabla_\xi\eta$ and
$\nabla_\xi f\eta=(\xi\cdot f)\eta+f\nabla_\xi\eta$ for $f\in
C^\infty(M,\Bbb R)$. This means that $\nabla$ defines a
\textit{partial connection} on $T^{-1}M$. 

Via the isomorphism $\La^2T^{-1}M\cong\gr_{-2}(TM)$, we obtain an
induced partial connection on the line bundle $\gr_{-2}(TM)$. This is
an operator $\nabla:\Ga(T^{-1}M)\times\Ga(\gr_{-2}(TM))\to
\Ga(\gr_{-2}(TM))$, which is linear over smooth functions in the first
variable and satisfies a Leibniz rule in the second variable. The
induced connection is characterized by
\begin{equation}
  \label{nabla-2-def}
  \nabla_{\gamma}\{\xi,\eta\}=\{\nabla_{\gamma}\xi,\eta\}+\{\xi,
\nabla_{\gamma}\eta\},
\end{equation}
for $\xi,\eta,\ga\in\Ga(T^{-1}M)$. 

This immediately singles out a preferred class of fields $r$ as above,
namely those, for which the nowhere vanishing section $\ph=q_{-2}(r)$
is parallel for the induced partial connection. We can completely
describe all such fields:
\begin{prop}\label{prop22}
  Let $\ph$ be a local non--vanishing smooth section of
  $\gr_{-2}(TM)$. Then there is a unique smooth section
  $r\in\Ga(T^{-2}M)$ such that $q_{-2}(r)=\ph$ and such that $\ph$ is
  parallel for the partial connection determined by $r$.
\end{prop}
\begin{pf}
  Since $\gr_{-2}(TM)$ is a quotient bundle of $T^{-2}M$ we find a
  local section $r_0\in\Ga(T^{-2}M)$ such that $q_{-2}(r_0)=\ph$. Any
  other section with this property is of the form $r_0+\delta$ for
  some $\delta\in\Ga(T^{-1}M)$. We have to compute how the choice of
  $\delta$ influences the partial connection. We will denote partial
  connections associated to $r_0$ by $\nabla$ and those associated to
  $r_0+\delta$ by $\nabla^\delta$. We first compute what happens on
  $T^{-1}M$. From the defining equation \eqref{nabla-1-def}, we
  immediately get
$$
\{\nabla^\delta_\ga\xi,\ph\}=\{\nabla_\ga\xi,\ph\}+\{\ga,\{\xi,\delta\}\}, 
$$
for all $\xi,\ga\in\Ga(T^{-1}M)$.  Now we can define a skew
symmetric bilinear map $a:T^{-1}M\x_M T^{-1}M\to\Bbb R$ by
$\{\xi,\ga\}=a(\xi,\ga)\ph$ for all $\xi,\ga\in\Ga(T^{-1}M)$. Then the
above equation shows that
$\nabla^\delta_\ga\xi=\nabla_\ga\xi+a(\xi,\delta)\ga$. 
Now choose local smooth sections $\xi,\eta\in\Ga(T^{-1}M)$ such
that $\ph=\{\xi,\eta\}$. Using the defining equation
\eqref{nabla-2-def} we immediately conclude that
$$
\nabla^\delta_\ga\{\xi,\eta\}=\nabla_\ga\{\xi,\eta\}+
a(\xi,\delta)\{\ga,\eta\}+a(\eta,\delta)\{\xi,\ga\}. 
$$
For fixed $\xi$, $\eta$, and $\delta$, the last two terms in the right
hand side define a bundle map $T^{-1}M\to\gr_{-2}(TM)$.  For
$\ga=\xi$, we obtain
$a(\xi,\delta)\{\xi,\eta\}=a(\xi,\delta)\ph=\{\xi,\delta\}$.
Likewise, inserting $\ga=\eta$, we obtain $\{\eta,\delta\}$. Since
$\{\xi,\eta\}=\ph$ is nowhere vanishing, the two fields form a local
frame for $T^{-1}M$. Hence the bundle map reduces to $\{\ga,\delta\}$
for any $\ga$, and locally
$\nabla^\delta_\ga\ph=\nabla_\ga\ph+\{\ga,\delta\}$, which then has to
hold globally. But $\ga\mapsto-\nabla_\ga\ph$ is a bundle map
$T^{-1}M\to\gr_{-2}(TM)$, so there is a unique section
$\delta\in\Ga(T^{-1}M)$ such that $-\nabla_\ga\ph=\{\ga,\delta\}$ for
all $\ga\in\Ga(T^{-1}M)$.\qed
\end{pf}

The proposition immediately implies that for a generalized contact
form $\al$ for $\Cal H$, there is a unique section $r\in\Ga(T^{-2}M)$
such that $\al(r)=1$ and such that $q_{-2}(r)$ is parallel for the
partial connection determined by $r$. This section is called the
\textit{generalized Reeb field} associated to $\al$. 

Notice that a generalized contact form $\al$ can be equivalently
viewed as a section of the bundle $\gr_2(T^*M)=(\gr_{-2}(TM))^*$. Any
partial connection on $\gr_{-2}(TM)$ induces a partial connection on
the dual bundle. Then a section $r\in\Ga(T^{-2}M)$ such that
$\al(r)=1$ is the generalized Reeb field if and only if
$\al\in\Ga(\gr_2(T^*M))$ is parallel for the partial connection
induced by $r$.

\subsection{The canonical extension of a generalized contact
  form}\label{2.3} 
Using the generalized Reeb field we can next construct a canonical
extension of a generalized contact form to a true one--form on $M$.

\begin{prop}\label{prop23}
  Let $\alpha$ be a generalized contact form with generalized Reeb
  field $r$. Then there is a unique one form $\tilde{\al}\in\Om^1(M)$
  extending $\al$ such that $i_rd\tilde{\alpha}|_{T^{-1}M}=0$.
\end{prop}
\begin{pf}
For a vector field $\zeta\in\frak X(M)$, we can write
$q_{-3}(\zeta)\in\Ga(\gr_{-3}(TM))$ as $\{r,\ze_1\}=q_{-3}([r,\ze_1])$
for a unique $\ze_1\in\Ga(T^{-1}M)$. But this exactly means that
$\ze-[r,\ze_1]\in\Ga(T^{-2}M)$. Hence for any extension $\tilde\al$ of
$\al$, we obtain
$$
\tilde\al(\ze)=\tilde\al([r,\ze_1])+\al(\ze-[r,\ze_1])=
-i_rd\tilde\al(\ze_1)+\al(\ze-[r,\ze_1]). \qed
$$
\end{pf}

In the sequel, we will use the same symbol $\al$ to denote a
generalized contact form and its canonical extension to a one form on
$M$.

\subsection{The projection associated to a generalized contact form}\label{2.4}
The objects we have associated to a generalized contact form $\al$ so
far amount to a partial splitting of the filtration of the tangent
bundle. Denoting by $r$ the generalized Reeb field associated to
$\al$, we obtain a projection from $T^{-2}M$ onto the subbundle
$T^{-1}M$ by $\ze\mapsto \ze-\al(\ze)r$. Likewise, the canonical
extension $\al\in\Om^1(M)$ of the generalized contact form allows one
to project an arbitrary tangent vector onto a multiple of $r$. To
construct a metric, we need a complete splitting of the filtration.
Hence we need a projection $\pi_{-1}:TM\to T^{-1}M$, which extends the
one on $T^{-2}M$ from above.

As we have noted in the proof of Proposition \ref{prop23}, given a vector
field $\ze\in\frak X(M)$, we find a unique $\ze_1\in\Ga(T^{-1}M)$ such
that $\ze-[r,\ze_1]\in\Ga(T^{-2}M)$. To decompose in a slightly finer
way, observe that this implies that
$\ze_2:=\ze-[r,\ze_1]-\al(\ze)r\in\Ga(T^{-1}M)$, so we can write
$$
\ze=[r,\ze_1]+\al(\ze)r+\ze_2
$$
for uniquely determined $\ze_1,\ze_2\in\Ga(T^{-1}M)$. In particular,
for any projection $\pi_{-1}$ as above, we get
$\pi_{-1}(\ze)=\pi_{-1}([r,\ze_1])+\ze_2$. Such a projection is
therefore equivalent to the operator $A:\Ga(T^{-1}M)\to\Ga(T^{-1}M)$
defined by $A(\xi):=\pi_{-1}([r,\xi])$. If we want $\pi_1$ to be
linear over smooth functions, we have to require that
$A(f\xi)=fA(\xi)+(r\cdot f)\xi$ for all $f\in C^\infty(M,\Bbb R)$.

There are two candidates for such an operator. First, define
$\Phi:\Gamma(T^{-1}M)\to\Gamma(T^{-1}M)$ by
\begin{align}\label{Phidef}
\nabla_{\xi}\nabla_{\eta}\gamma-\nabla_{\eta}\nabla_{\xi}\gamma-
\nabla_{[\xi,\eta]-\al([\xi,\eta])r}\gamma=\alpha([\xi,\eta])\Phi(\gamma),
\end{align}
for $\xi,\eta,\ga\in\Ga(T^{-1}M)$. This makes sense since the left
hand side of \eqref{Phidef} is alternating in $\xi$ and $\eta$ and
hence depends only on $\alpha([\xi,\eta])$. Note further that
$[\xi,\eta]-\al([\xi,\eta])r$ is the projection of
$[\xi,\eta]\in\Ga(T^{-2}M)$ to a section of $T^{-1}M$. Extending
$\nabla$ to a linear connection on $TM$ and denoting by $R$ the
curvature of this extension, the left hand side of \eqref{Phidef}
(which is evidently independent of the extension) can be written as
$R(\xi,\eta)(\ga)+\al([\xi,\eta])\nabla_r\ga$. This shows that
$\Phi(f\gamma)=f\Phi(\gamma)+(r\cdot f)\gamma$.

\noindent
Second, $\ga\mapsto q_{-3}([r,[\gamma,r]])$ defines an operator
$\Ga(T^{-1}M)\to\Ga(T^{-3}M)$. Since $q_{-2}(r)$ is nowhere vanishing,
there is a unique operator
$\Psi:\Gamma(T^{-1}M)\to\Gamma(T^{-1}M)$ such that
\begin{align}\label{Psidef}
\{\Psi(\gamma),q_{-2}(r)\}=\tfrac{1}{2}q_{-3}([r,[\gamma,r]]).
\end{align}
From this definition, $\Psi(f\gamma)=f\Psi(\gamma)+(r\cdot f)\gamma$
follows easily. 

Any convex combination of the operators $\Ph$ and $\Ps$ also has the
right behavior under multiplication by smooth functions. To define our
projection $\pi_{-1}:TM\to T^{-1}M$ we use the combination
$\tfrac{-2}{5}\Ph+\tfrac{7}{5}\Ps$, so we define
\begin{align}\label{pidef}
  \pi_{-1}(\zeta):=\tfrac{-2}{5}\Phi(\zeta_1)+
\tfrac{7}{5}\Psi(\zeta_1)+\zeta_2,
\end{align} 
where $\zeta_1$ is the unique vector field in $\Gamma(T^{-1}M)$ such
that $q_{-3}(\zeta)=q_{-3}([r,\zeta_1])$ and
$\zeta_2=\zeta-[r,\zeta_1]-\alpha(\zeta) r$. The motivation for the
choice of factors will become clear from the following computations
and from section \ref{3}.

We can now define the pseudo--Riemannian metric associated to $\al$.
Let $\zeta,\zeta'$ be tangent vectors on $M$. Using the components
$\ze_1$ and $\ze_2$ as above and likewise for $\ze'$, we define
\begin{equation}\label{metric}
  g_{\alpha}(\ze,\ze'):=d\alpha(\ze_1,\pi_{-1}(\ze'))-
\tfrac{4}{3}\al(\ze)\al(\ze')+d\al(\ze_1',\pi_{-1}(\ze)).
\end{equation}
\begin{prop}\label{prop24}
  For any generalized contact form $\al$, the map $g_\al$ defined
  above is a pseudo--Riemannian metric of signature $(2,3)$ on $M$. 
\end{prop}
\begin{pf}
  Evidently, $g_{\al}$ is a smooth, symmetric bilinear bundle map.  For
  $\ze\in T^{-1}M$ we have $\ze_1=\al(\ze)=0$ and $\pi_{-1}(\ze)=\ze$.
  This shows that the rank two subbundle $T^{-1}M\subset TM$ is
  isotropic for $g_\al$. Likewise, $\ker(\pi_{-1})\cap\ker(\al)$ is a
  rank two subbundle of $TM$, which is transversal to $T^{-1}M$ and
  isotropic for $g_{\al}$. Taking $\ze'$ in
  $\ker(\pi_{-1})\cap\ker(\al)$ and $\ze\in T^{-1}M$, we by definition
  get $g_{\alpha}(\ze,\ze')=d\al(\ze_1',\ze)$. This vanishes for all
  $\ze$ if and only if $\ze_1'=0$ and hence $\ze'\in T^{-2}M$. But
  then $\al(\ze')=0$ implies $\ze'\in T^{-1}M$ and $\pi_{-1}(\ze')=0$
  shows $\ze'=0$. Hence $g_\al$ induces a non--degenerate pairing
  between the two isotropic subbundles. Since $r$ spans a rank one
  subbundle transversal to the two isotropic subbundles on which
  $g_\al$ is negative definite, the result follows.
\qed\end{pf}

\subsection{The dependence on the generalized contact form}\label{2.5}
To analyze how $g_\al$ depends on $\al$, we have to study the dependence
of the ingredients used in the construction. Given a generalized
contact form $\al$, any other generalized contact form is obtained by
multiplying $\al$ by a nowhere vanishing smooth function. We will
following the convention that we denote the changed generalized
contact form by $\hat\al$ and indicate all quantities referring to the
new form by a hat. 

Let us first check what happens if we replace $\al\in\Ga((T^{-2}M)^*)$
by $\hat\al:=-\al$. If in the defining equation \eqref{nabla-1-def} we
replace $r$ by $-r$ (and hence $\ph$ by $-\ph$) we obtain the same
partial connection $\nabla$. This shows that $\hat r=-r$ and
$\hat\nabla=\nabla$. Then it follows that $\hat\al=-\al$ also holds
for the extensions to one forms on $M$. Decomposing $\ze\in\frak X(M)$
as introduced in \ref{2.4}, we obtain $\hat\ze_1=-\ze_1$ and
$\hat\ze_2=\ze_2$. In the defining equation \eqref{Phidef}, the left
hand side remains unchanged while in the right hand side $\al$ has to
be replaced by $-\al$, so $\hat\Ph=-\Ph$. Similarly, the definition in
\eqref{Psidef} shows that $\hat\Ps=-\Ps$. Hence we obtain
$\hat\pi_{-1}=\pi_{-1}$. Putting all these results together, we
conclude that $g_{-\al}=g_{\al}$ from the definition of the metric.

Hence it suffices to analyze the behavior of $g_{\al}$ under rescaling
$\al$ by a positive smooth function, which we write as $e^f$ for $f\in
C^\infty(M,\Bbb R)$.

\begin{lem}\label{lem25}
  Let $\al\in\Ga((T^{-2}M)^*)$ be a generalized contact form, consider
  a smooth function $f\in C^\infty (M,\Bbb R)$, and the generalized
  contact form $\hat\al= e^f\al\in\Ga((T^{-2}M)^*)$. Then we have:

\noindent
(i) The Reeb vector field $\hat r$ associated to $\hat\al$ is given by
$\hat{r}=e^{-f}r+\delta$, where $\delta\in\Gamma(T^{-1}M)$ is the
unique vector field such that $\{\ga,\delta\}=4e^{-f}df(\ga)q_{-2}(r)$
for all $\gamma\in\Gamma(T^{-1}M)$.

\noindent
(ii) The canonical extensions of the two generalized contact forms to
one forms on $M$ are related by
$\hat{\al}(\ze)=e^f\al(\ze)+3df(\zeta_1)$, where
$\ze_1\in\Ga(T^{-1}M)$ is characterized by
$q_{-3}(\ze)=\{\ph,\ze_1\}.$
\end{lem}
 
\begin{pf}
  (i) Put $r_0=e^{-f}r$. Then $\hat\al(r_0)=1$, and we can compute
  $\hat r$ following the proof of Proposition \ref{prop22}. We first have
  to compute the partial connection $\tilde\nabla$ induced by $r_0$.
  Putting $\ph=q_{-2}(r)$ and using that for $\ga\in\Ga(T^{-1}M)$ we
  have $\ga\cdot e^{-f}=-e^{-f}df(\ga)$ one easily computes that
$$
q_{-3}([\ga,[\xi, r_0]])=e^{-f}\bigg(q_{-3}([\ga,[\xi,r]])-
df(\ga)\{\xi,\ph\}-df(\xi)\{\ga,\ph\}\bigg). 
$$
Via the defining equation \eqref{nabla-1-def} in \ref{2.2} and
$q_{-2}(r_0)=e^{-f}\ph$, this shows that the partial connection
$\tilde\nabla$ on $T^{-1}M$ determined by $r_0$ is given by 
$$
\tilde\nabla_\ga\xi=\nabla_\ga\xi-df(\ga)\xi-df(\xi)\ga. 
$$
For the induced connection on $\gr_{-2}(TM)$ we therefore get
$$
\tilde\nabla_\ga\{\xi,\eta\}=\nabla_\ga\{\xi,\eta\}-
2df(\ga)\{\xi,\eta\}-df(\xi)\{\ga,\eta\}-df(\eta)\{\xi,\ga\}. 
$$
Fixing $\xi$ and $\eta$, the last two terms in the right hand side
define a bundle map $T^{-1}M\to\gr_{-2}(TM)$. For $\ga=\xi$ and
$\ga=\eta$, this bundle map coincides with $-df(\ga)\{\xi,\eta\}$. To
have $\{\xi,\eta\}\neq 0$, the two sections have to form a frame of
$T^{-1}M$, so we conclude that the induced connection on
$\gr_{-2}(TM)$ is characterized by
$\tilde\nabla_\ga\ph=\nabla_\ga\ph-3df(\ga)\ph$.  Using
$\nabla_\ga\ph=0$, this implies
$$
\tilde\nabla_\ga e^{-f}\ph=-4df(\ga)e^{-f}\ph,
$$
and (i) follows from the proof of Proposition \ref{prop22}.

\noindent
(ii) From (i) we know that $\hat r=e^{-f}r+\delta$.  For
$\ga\in\Ga(T^{-1}M)$, we thus have
\begin{equation}
  \label{klammer}
[\hat r,\ga]=e^{-f}[r,\ga]+e^{-f}df(\ga)r+[\delta,\ga].
\end{equation}
Since $\al([r,\ga])=0$, the one form $e^f\al\in\Om^1(M)$ maps this to
$df(\ga)-e^f\al([\ga,\delta])$. By definition,
$\{\ga,\delta\}=4e^{-f}df(\ga)\ph$, and hence
$\al([\ga,\delta])=4e^{-f}df(\ga)$. Thus we obtain
$$
e^f\al([\hat r,\ga])=-3df(\ga).
$$
But this exactly means that the claimed formula for $\hat\al$ defines
a form which annihilates each field of the form $[\hat r,\ga]$ for
$\ga\in\Ga(T^{-1}M)$. Since it obviously is an extension of
$\hat\al\in\Ga((T^{-2}M)^*)$,this completes the proof.
\qed\end{pf}
\begin{rem}\label{rem25}
  Note that the transformation laws in the lemma both depend only on
  $df|_{T^{-1}M}$, i.e.~on the class of $df$ in $\Ga(\gr_1(T^*M))$.
\end{rem}

It remains to analyze the dependence of $\pi_{-1}$ on the generalized
contact form. 

\begin{lem}\label{lem26}
  Let $\al\in \Ga((T^{-2}M)^*)$ be a generalized contact form and
  consider a rescaling $\hat\al=e^f\al$ for $f\in C^\infty(M,\Bbb R)$.
  Let $r$ be the generalized Reeb field for $\al$ and put
  $\ph=q_{-2}(r)$. Then the projection $\hat\pi_1$ determined by
  $\hat\al$ is given by
  $$
  \hat\pi_{-1}(\ze)=\pi_{-1}(\zeta)-\tfrac{3}{2}df(r)\ze_1-\tfrac{3}{2}e^fdf(\ze_1)\delta-e^f\al(\ze)\delta,
  $$
  where $\ze_1\in\Ga(T^{-1}M)$ is the unique vector field such that
  $q_{-3}(\ze)=\{\ph,\ze_1\}$, and $\delta\in\Ga(T^{-1}M)$ is
  characterized by $\{\ga,\delta\}=4e^{-f}df(\ga)q_{-2}(r)$ for all
  $\gamma\in\Gamma(T^{-1}M)$.
\end{lem}
\begin{pf}
  Let us first compare the decompositions of $\ze\in\frak X(M)$ from
  \ref{2.4} with respect to the two generalized contact forms. For
  $\al$, this reads as 
$$
\ze=[r,\ze_1]+\al(\ze)r+\ze_2,
$$
where $\ze_1\in\Ga(T^{-1}M)$ is characterized by
$q_{-3}(\ze)=\{\ph,\ze_1\}$, and then $\ze_2\in T^{-1}M$ is defined by
the equation. Since $\hat\ph=e^{-f}\ph$ we see that
$\widehat{\ze_1}=e^f\ze_1$. From Lemma \ref{lem25}, we know that $\hat
r=e^{-f}r+\delta$. Note that by definition of $\delta$, we get
$df(\delta)=0$. Using this and formula \eqref{klammer} from the proof
of Lemma \ref{lem25}, we obtain
\begin{equation}
  \label{r-zeta}
  [\hat r,\hat\ze_1]=[r,\ze_1]+df(r)\ze_1+df(\ze_1)r+e^f[\delta,\ze_1]. 
\end{equation}
Further, the formulae from Lemma \ref{lem25} show that
$$
\hat\al(\ze)\hat r=(\al(\ze)+3df(\ze_1))(r+e^f\delta).
$$
Putting the results obtained so far together, we get
\begin{equation}
  \label{hatze2}
  \widehat{\ze_2}-\ze_2=-e^{f}\al(\ze)\delta-df(r)\ze_1-
3e^{f}df(\ze_1)\delta-4df(\ze_1)r- e^f[\delta,\zeta_1].
\end{equation}
Using \eqref{r-zeta}, we next compute 
\begin{align*}
 q_{-3}([\hat{r},[\hat{\zeta}_1,\hat{r}]])
 =&e^{-f}q_{-3}([r,[\ze_1,r]])+2q_{-3}([\delta,[\zeta_1,r]])-
q_{-3}([\ze_1,[\delta,r]])\\
-&e^f\{\delta,\{\delta,\ze_1\}\}-e^{-f}df(r)\{\ph,\ze_1\}-
df(\ze_1)\{\delta,\ph\}.
\end{align*}
The second and third term can be rewritten in terms of the partial
connection $\nabla$ associated to $r$ using the definition in formula
\eqref{nabla-1-def} in \ref{2.2}. On the other hand,
$\{\delta,\ze_1\}=-4e^{-f}df(\ze_1)\ph$ by definition. In view of
equation \eqref{Psidef} from \ref{2.4} this shows that
\begin{equation}
  \label{hatPsi}
  \hat\Psi(\hat\ze_1)-\Psi(\ze_1)=e^f\nabla_\delta\ze_1-
\tfrac{1}{2}e^f\nabla_{\ze_1}\delta+\tfrac{3}{2}e^fdf(\ze_1)\delta+
\tfrac{1}{2}df(r)\ze_1.
\end{equation}
To compute $\hat\Ph(\hat\ze_1)$ we have to analyze the relation
between the partial connections associated to $\al$ and $\hat\al$. 
Using \eqref{klammer} from the proof of Lemma \ref{lem25}, one computes
$q_{-3}([\xi,[\eta,\hat r]])$, and via the defining equation
\eqref{nabla-1-def} from \ref{2.2} this shows that 
$$
\hat\nabla_\xi\eta=\nabla_\xi\eta-df(\xi)\eta+3df(\eta)\xi.
$$
Using this, one easily verifies directly that for
$\xi,\eta,\ga\in\Ga(T^{-1}M)$, the difference
$\hat\nabla_\xi\hat\nabla_\eta\ga-\nabla_\xi\nabla_\eta\ga$ can, up to
terms symmetric in $\xi$ and $\eta$, be expressed as
\begin{equation}
  \label{double-deriv}
  \begin{aligned}
      3df(\nabla_\eta\ga)\xi&-(\xi\cdot df(\eta))\ga+3(\xi\cdot
  df(\ga))\eta\\
&+3df(\ga)\nabla_\xi\eta+9df(\ga)df(\eta)\xi
  \end{aligned}
\end{equation}
Now $[\xi,\eta]\in\Ga(T^{-2}M)$, which implies that 
$$
\hat\al([\xi,\eta])\hat r=\al([\xi,\eta])(r+e^f\delta).
$$ 
Using this and $df(\delta)=0$, one shows that 
$\hat\nabla_{[\xi,\eta]-\hat\al([\xi,\eta])\hat
  r}\ga-\nabla_{[\xi,\eta]-\al([\xi,\eta])r}\ga$ is given by
\begin{equation}
  \label{bracket-deriv}
  \begin{aligned}
  &-df([\xi,\eta]-\al([\xi,\eta])r)\ga+3df(\ga)([\xi,\eta]-
\al([\xi,\eta])r)\\
&-\hat\al([\xi,\eta])\left(\nabla_\delta\ga+3df(\ga)\delta\right).     
  \end{aligned}
\end{equation}
Now we need a few identities. First, expanding $0=ddf(\xi,\eta)$, we
obtain 
\begin{equation}
  \label{d2}
  \xi\cdot df(\eta)-\eta\cdot
  df(\xi)-df([\xi,\eta]-\al([\xi,\eta])r)=\al([\xi,\eta])df(r). 
\end{equation}
For $\ga_1,\ga_2\in\Ga(T^{-1}M)$, the Jacobi identity implies
\begin{align*}
q_{-3}([\ga_1,[\ga_2,r]])-q_{-3}([\ga_2,[\ga_1,r]])=q_{-3}([[\ga_1,\ga_2],r]).
\end{align*}
The right hand side remains unchanged if we replace $[\ga_1,\ga_2]$ by
$[\ga_1,\ga_2]-\al([\ga_1,\ga_2])r\in\Ga(T^{-1}M)$ and then gives
$\{[\ga_1,\ga_2]-\al([\ga_1,\ga_2])r,\phi\}$. Hence
\begin{equation}\label{torsfree}
\nabla_{\gamma_1}\gamma_2-\nabla_{\gamma_2}\gamma_1=
[\gamma_1,\gamma_2]-\alpha([\gamma_1,\gamma_2])r,
\end{equation}
which is the analog of torsion freeness for $\nabla$. To obtain a
formula for $\hat\Phi(\ga)-\Ph(\ga)$, we have to take
\eqref{double-deriv}, then subtract the analogous terms with $\xi$ and
$\eta$ exchanged and further subtract \eqref{bracket-deriv}. Using
\eqref{d2} and \eqref{torsfree}, we obtain
\begin{equation}
  \label{hatPhi1}
  \begin{aligned}
    3df(\nabla_\eta\ga)\xi-&3df(\nabla_\xi\ga)\eta+
3(\xi\cdot df(\ga))\eta-3(\eta\cdot df(\ga))\xi\\
-&\al([\xi,\eta])df(r)\ga+9df(\ga)(df(\eta)\xi-df(\xi)\eta)\\
    +&\hat\al([\xi,\eta])(\nabla_\delta\ga+3df(\ga)\delta).
  \end{aligned}
\end{equation}
Inserting \eqref{torsfree} into \eqref{d2}, we get
$$
df(\nabla_\eta\ga)-\eta\cdot df(\ga)=df(\nabla_\ga\eta)-\ga\cdot
df(\eta)-\al([\eta,\ga])df(r). 
$$ 
Using this and the analogous formula for $\xi$ and $\ga$, we see that
the first line of \eqref{hatPhi1} can be rewritten as
\begin{equation}
  \label{hatPhi2}
  \begin{aligned}
    3(\ga\cdot df(\xi))\eta-&3df(\nabla_\ga\xi)\eta-
3(\ga\cdot df(\eta))\xi+3df(\nabla_\ga\eta)\xi\\
+&3df(r)(\al([\xi,\ga])\eta-\al([\eta,\ga])\xi).
  \end{aligned}
\end{equation}
Now $\ga\mapsto \al([\xi,\ga])\eta-\al([\eta,\ga])\xi$ defines a
bundle map $T^{-1}M\to T^{-1}M$. This map coincides with
$\al([\xi,\eta])\ga$ for $\ga=\xi$ and $\ga=\eta$, and since
$\{\xi,\eta\}$ is nowhere vanishing by assumption, this holds for all
$\ga$. 

Similarly, we see that the bundle map $T^{-1}M\to\mathrm{gr}_{-2}(TM)$
given by $\gamma\mapsto df(\xi)\{\gamma,\eta\}-df(\eta)\{\gamma,\xi\}$
coincides with $df(\gamma)\{\xi,\eta\}$. By definition of $\delta$,
this implies
\begin{equation}
   \label{d1f}
   4(df(\xi)\eta-df(\eta)\xi)=
 \al([\xi,\eta])e^f\delta. 
 \end{equation}
It follows that 
\begin{align*}
 4(df(\nabla_\ga\xi)&\eta-
 df(\eta)\nabla_\ga\xi+df(\xi)\nabla_\ga\eta-df(\nabla_\ga\eta)\xi)=\\
&\alpha([\nabla_{\gamma}\xi,\eta])e^f\delta-\alpha([\nabla_{\gamma}\eta,\xi])e^f\delta.
\end{align*}
Next we use the characterisation \eqref{nabla-2-def} of the induced
connection on $\mathrm{gr}_{-2}(TM)$, the fact that
$\nabla_{\gamma}\phi=0$ and
$\{\gamma_1,\gamma_2\}=\alpha([\gamma_1,\gamma_2])\phi$ for all
$\gamma_1,\gamma_2\in\Gamma(T^{-1}M)$ to get
$$
\alpha([\nabla_{\gamma}\xi,\eta])e^f\delta-
\alpha([\nabla_{\gamma}\eta,\xi])e^f\delta=
(\gamma\cdot\alpha([\xi,\eta]))e^f\delta.
$$
Now we apply $\nabla_\ga$ to \eqref{d1f} and simplify using the last
two equations to obtain
\begin{equation}
  \label{helpful}
  \begin{aligned}
df(\nabla_\ga\eta)\xi&-(\ga\cdot
df(\eta))\xi-df(\nabla_\ga\xi)\eta+(\ga\cdot
df(\xi))\eta\\
&=\tfrac{1}{4}\hat\al([\xi,\eta])(df(\ga)\delta+
\nabla_\ga\delta). 
  \end{aligned}
\end{equation}
Collecting our results and using the defining equation \eqref{Phidef}
from \ref{2.4}, we see that 
$$
\hat\Ph(\ga)=\Ph(\ga)+2e^{-f}df(r)\ga+\tfrac{3}{2}df(\ga)\delta+
\nabla_\delta\ga+\tfrac{3}{4}\nabla_\ga\delta.  
$$
Inserting $\ga=\hat\ze_1=e^f\ze_1$, and using $df(\delta)=0$ we get
\begin{equation}
  \label{hatPhi3}
\hat\Ph(\hat\ze_1)-\Ph(\ze_1)=3df(r)\ze_1+\tfrac{3}{2}e^fdf(\ze_1)\delta+e^f\nabla_\delta\ze_1+
\tfrac34e^f\nabla_{\ze_1}\delta 
\end{equation}
Using appropriate multiples of this, \eqref{hatPsi}, and
\eqref{hatze2} we obtain the following expression for
$\hat\pi_1(\ze)-\pi_1(\ze)$: 
\begin{equation}\label{final}
  e^f(\nabla_\delta\ze_1-\nabla_{\ze_1}\delta-[\delta,\ze_1])-4df(\ze_1)r
-\tfrac{3}{2}df(r)\ze_1-\tfrac{3}{2}e^fdf(\ze_1)\delta-e^f\al(\ze)\delta. 
\end{equation}
Using \eqref{torsfree}, the first four terms simplify to 
$$
-(e^f\al([\delta,\ze_1])+4df(\ze_1))r,
$$
and we have seen in the proof of Lemma \ref{lem25} that this
vanishes. \qed
\end{pf}

\subsection{The main result}
Having the technical results at hand, it is now easy to prove that a
change of generalized contact form just leads to a conformal rescaling
of the associated pseudo--Riemannian metric. 

\begin{thm}\label{thm26}
  Replacing a generalized contact form $\alpha\in\Ga(T^{-2}M)$ by
  $\hat{\alpha}=e^f\alpha$, the pseudo--Riemannian metrics associated
  to $\al$ and $\hat\al$ as in \eqref{metric} in \ref{2.4} are related
  by $g_{\hat{\alpha}}=e^{2f}g_{\alpha}$. In particular, the conformal
  class of $g_\al$ depends only on the generic distribution $\Cal H$.
\end{thm}
\begin{pf}
  In the definition of $g_\al$, the exterior derivative $d\al$ is only
  applied to two elements of $T^{-1}M$, whence passing to $\hat\al$
  this only rescales by $e^f$. Consequently, the Lemma shows that
  $d\hat\al(\hat\ze_1,\hat\pi_1(\ze'))-e^{2f}d\al(\ze_1,\pi_1(\ze'))$
  is given by
$$
-e^{2f}\bigg(\tfrac{3}{2}df(r)d\al(\ze_1,\ze_1')+
e^f(\al(\ze')+\tfrac{3}{2}df(\ze_1'))d\al(\ze_1,\delta)\bigg).
$$
The first summand in the bracket is skew symmetric in $\ze$ and $\ze'$
and hence will not contribute to the final result. In the end of the
proof of the Lemma, we have seen that 
$$
e^fd\al(\ze_1,\delta)=-e^f\al([\ze_1,\delta])=-4df(\ze_1). 
$$
Inserting this and using part (ii) of Lemma \ref{lem25}, the result
follows by a simple direct computation.\qed
\end{pf}

\section{The relation to Nurowski's construction}\label{3}
In this section, we will describe the relation of the conformal class
constructed in section \ref{2} to the canonical Cartan connection
associated to a generic rank two distribution on a
five--manifold. In particular, this will show that our conformal class
coincides with the one constructed by P. Nurowski in
\cite{Nurowski}. Moreover, this will put our construction in a
broader context of general tools for parabolic geometries which have
been developed during the last years. We will start by describing the
canonical Cartan connection associated to a generic rank two
distribution in dimension five. 

\subsection{On $G_2$}\label{3.1}
By Cartan's classical result \cite{Cartan:five}, generic rank two
distributions in dimension five admit a canonical Cartan connection on
a certain principal bundle. The structure group of this bundle is a
subgroup in a Lie group whose Lie algebra is the split real form of
the exceptional Lie algebra of type $G_2$. We next discuss the
necessary background on this Lie algebra put Cartan's result into the
perspective of the general theory of parabolic geometries.

A Lie group $G$ with this Lie algebra can be realized as the
automorphism group of the algebra $\Bbb O_s$ of split octonions, see
\cite{Springer}. $\Bbb O_s$ is an eight dimensional non--associative
unital real algebra with a multiplicative inner product of split
signature $(4,4)$. Any automorphism of $\Bbb O_s$ preserves the inner
product and the unit element, so the automorphism group naturally acts
on the orthocomplement of the unit element. This is the seven
dimensional space $\im(\Bbb O_s)$ of purely imaginary split
octonions, which carries an invariant inner product of signature
$(3,4)$. Hence $G$ can be naturally viewed as a closed subgroup of
$SO(3,4)$. The Lie algebra $\frak g$ of $G$ is the algebra of
derivations of $\Bbb O_s$. Since any derivation vanishes on the unit
element, also $\fg$ is naturally represented on $\im(\Bbb O_s)$ and
$\fg\subset\frak{so}(3,4)$.

To obtain an explicit description of $\frak g$, we first fix the inner
product of signature $(3,4)$. In terms of coordinates $x_0,\dots,x_6$
on $\Bbb R^7$ consider the quadratic form
$x_0x_6+x_1x_4+x_2x_5-(x_3)^2$, which is evidently induced by an inner
product of signature $(3,4)$. The explicit form of $\frak g$ for this
inner product can be essentially read off from \cite{Sa1}. In the
notation of that article, one has to use the ordered basis
$\{X_1,X_6,X_7,X_4,X_2,X_3,X_5\}$ to obtain
$$\mathfrak{g}=\left\{\begin{pmatrix}
    \mathrm{tr}(A)&Z&s&W&0\\
    X&A&\sqrt{2}\mathbb{J}Z^t&\frac{s}{\sqrt{2}}\mathbb{J}&-W^t\\
    r&-\sqrt{2}X^t\mathbb{J}&0&-\sqrt{2}Z\mathbb{J}&s\\
    Y&-\frac{r}{\sqrt{2}}\mathbb{J}&\sqrt{2}\mathbb{J}X&-A^t&-Z^t\\
    0&-Y^t&r&-X^t&-\mathrm{tr}(A)
\end{pmatrix}\right\}$$
with $A\in\mathfrak{gl}(2,\mathbb{R})$, $X,Y,Z^t,W^t\in\mathbb{R}^2$, 
$r,s\in\mathbb{R}$ and $\mathbb{J}:=\left(\begin{smallmatrix}0&-1\\
1&0\end{smallmatrix}\right)$. Indeed, one may easily verify directly
that this forms a Lie subalgebra of $\frak{so}(3,4)$, the diagonal
matrices contained in $\frak g$ act diagonalizably under the adjoint
action, and the resulting root decomposition of $\frak g$ has a root
system of type $G_2$. 

Let us decompose
$\mathfrak{g}=\mathfrak{g}_{-3}\oplus\mathfrak{g}_{-2}\oplus
\mathfrak{g}_{-1}\oplus\mathfrak{g}_0\oplus\mathfrak{g}_{1}\oplus
\mathfrak{g}_{2}\oplus\mathfrak{g}_{3}$
as in
$$\begin{pmatrix}
  \mathfrak{g}_{0}&\mathfrak{g}_{1}&\mathfrak{g}_{2}&\mathfrak{g}_{3}&0\\
  \mathfrak{g}_{-1}&\mathfrak{g}_{0}&\mathfrak{g}_{1}&
   \mathfrak{g}_{2}&\mathfrak{g}_{3}\\
  \mathfrak{g}_{-2}&\mathfrak{g}_{-1}&0&\mathfrak{g}_{1}&\mathfrak{g}_{2}\\
  \mathfrak{g}_{-3}&\mathfrak{g}_{-2}&\mathfrak{g}_{-1}&\mathfrak{g}_{0}&
   \mathfrak{g}_{1}\\
  0&\mathfrak{g}_{-3}&\mathfrak{g}_{-2}&\mathfrak{g}_{-1}&\mathfrak{g}_{0}
\end{pmatrix}.
$$
Then this is immediately seen to define a grading on $\fg$,
i.e.~$[\fg_i,\fg_j]\subset\fg_{i+j}$, where we agree that
$\fg_\ell=\{0\}$ for $|\ell|>3$. In particular, the Lie bracket
defines a representation of the subalgebra $\fg_0$ on each $\fg_i$,
which is compatible with the Lie brackets. It is easy to see that the
representation of $\fg_0$ on $\fg_{-1}\cong\Bbb R^2$ is faithful, so
we can use this representation to identify $\fg_0$ with
$\frak{gl}(2,\Bbb R)$. The Lie bracket induces isomorphisms
$\La^2\fg_{-1}\to\fg_{-2}$ and $\fg_{-1}\otimes\fg_{-2}\to\fg_{-3}$ of
$\fg_0$--modules, which already indicates the relation to generic rank
two distributions. On the other hand, $\fg_i$ is dual to $\fg_{-i}$ as
a $\fg_0$--module for $i=1,2,3$. A convenient way to express these
dualities is via mapping two matrices to $\frac16$ times the trace of
their product. Let us denote all these pairings by $B$. In the
notation introduced above, this is explicitly given by $(X,Z)\mapsto
ZX$, $(r,s)\mapsto\frac{1}{2}rs$, and $(Y,W)\mapsto \tfrac13 WY$.

For $i=-3,\dots,3$, we define $\fg^i:=\fg_i\oplus\dots\oplus\fg_3$.
This makes $\fg$ into a filtered Lie algebra,
i.e.~$[\fg^i,\fg^j]\subset\fg^{i+j}$, where we agree that
$\fg^\ell=\fg$ for $\ell<-3$ and $\fg^\ell=\{0\}$ for $\ell>3$. In
particular, $\frak p:=\fg^0$ is a parabolic subalgebra of $\frak g$
and $\frak p_+:=\fg^1$ is a nilpotent ideal in $\frak p$.

\subsection{The canonical Cartan connection}\label{3.2}
As we have seen above, we may view $G=\Aut(\Bbb O_s)$ as a closed
subgroup of $SO(3,4)$. Define $P\subset G$ to be the intersection of
$G$ with the stabilizer of the isotropic line spanned by the first
basis vector. By construction, $P$ corresponds to the Lie subalgebra
$\frak p\subset\fg$. For $g\in P$, the adjoint action $\Ad(g)$
preserves the filtration on $\fg$, i.e.~$\Ad(g)(\fg^i)\subset\fg^i$
for all $i$. Define $G_0\subset P$ as the subgroup of those $g$ for
which $\Ad(g)(\fg_i)\subset\fg_i$ for all $i$. Then the Lie algebra of
$G_0$ is $\fg_0$. On the other hand, it turns out that the exponential
map defines a diffeomorphism from $\frak p_+$ onto a closed normal
subgroup $P_+\subset P$, and $P/P_+$ is naturally isomorphic to
$G_0$. 

It is well known (and easy to see from the explicit form of $\fg$
above) that the 7--dimensional representation of $\fg$ defined by the
above matrix form is irreducible. By Schur's lemma, this implies that
the center of $G$ consists only of multiples of the identity matrix,
so since $G\subset SO(3,4)$, we see that $G$ has trivial center. Thus
the adjoint representation maps $G$ injectively into the group of
automorphisms of $\fg$ and thus $G_0$ is mapped injectively into the
group $\Aut_{\mathrm{gr}}(\fg)$ of automorphisms of the
\textit{graded} Lie algebra $\fg$. From \cite{Sa1} we see that the
latter group coincides with $\Aut_{\mathrm{gr}}(\fg_-)\cong
GL(\fg_{-1})$. 

On the other hand, one easily verifies directly that any invertible
linear map on $\fg_{-1}$ can be obtained from the adjoint action of
an element $g\in G_0$. Hence we conclude that $G_0\cong
GL(\fg_{-1})\cong \Aut_{\mathrm{gr}}(\fg_-)$. As explained in
\cite{Sa1}, this shows that parabolic geometries of type $(G,P)$ are
equivalent to filtrations of the tangent bundle such that the bundle
of symbol algebras is locally trivial and modelled on $\fg_-$, and
hence to generic rank two distributions in dimension five. 

This can be easily made more explicit. Suppose that $M$ is a five
dimensional smooth manifold endowed with a generic distribution $\Cal
H\subset TM$ of rank two, and consider the corresponding filtration
$\{T^iM\}$ of the tangent bundle as introduced in \ref{2.1}. Let
$p_0:\Cal G_0\to M$ be the linear frame bundle of $\Cal H=T^{-1}M$.
This has structure group $GL(2,\Bbb R)$ which we view as $G_0$. Since
all the components $\gr_i(TM)$ of the associated graded can be
constructed from $T^{-1}M$, they are naturally associated to $\Cal
G_0$. This can be expressed via $G_0$--equivariant partially defined
one--forms with values in $\fg_-$, which define a \textit{regular
  infinitesimal flag structure} in the sense of \cite[2.6]{Weyl}.

The prolongation procedures for parabolic geometries show that $\Cal
G_0$ naturally extends to a principal $P$--bundle $p:\Cal G\to M$,
which can be endowed with a canonical normal Cartan connection
$\om\in\Om^1(\Cal G,\frak g)$. The bundle $\Cal G_0$ can be recovered
from $\Cal G$ as $\Cal G/P_+$. Further, the generic distribution $\Cal
H$ and, more generally, the filtration of $TM$ can be recovered from
$(p:\Cal G\to M,\om)$. For a tangent vector $\xi\in T_xM$ choose any
$u\in\Cal G$ and $\tilde\xi\in T_u\Cal G$ such that
$T_up\cdot\tilde\xi=\xi$.  Then $\xi\in T^i_xM$ for $i=-1,-2$ if and
only if $\om(\tilde\xi)\in\fg^i$.  This is independent of the choices
by the defining properties of a Cartan connection an the fact that
each $\fg^i$ is a $P$--invariant subspace of $\fg$.

Normality is a condition on the curvature of the Cartan connection
$\om$. The detailed form of the condition is not important for our
purposes. The curvature of $\om$ is most easily viewed as the two form
$\Cal K\in\Om^2(\Cal G,\frak g)$ defined by 
$$
\Cal K(\xi,\eta):=d\om(\xi,\eta)+[\om(\xi),\om(\eta)]. 
$$
The only fact about the normality condition we will need in the sequel
is a restriction on the homogeneity of $\Cal K$. Namely, if
$Tp\cdot\xi\in T^iM$ and $Tp\cdot\eta\in T^jM$ for some
$i,j=-3,\dots,-1$, then $\Cal K(\xi,\eta)\in \fg^{i+j+4}$ (and for the
sequel even $i+j+3$ would be sufficient).
 
\subsection{Nurowski's conformal structure}\label{3.3}
We have realized the group $G$ as a subgroup of $\tilde G:=SO(3,4)$.
Denoting by $\tilde P\subset\tilde G$ the stabilizer of the isotropic
line generated by the first basis vector, we see that $P=G\cap\tilde
P$. Now it is well known that $\tilde G/\tilde P$ is the M\"obius
space $S^{2,3}$ with $\tilde G$ acting as the group of all conformal
isometries of the canonical (locally conformally flat) conformal
structure. Denoting by $\tg$ and $\tp$ the Lie algebras, this
conformal structure is induced by a conformal class of inner products
on $\tg/\tp$ which is invariant under the natural action of $\tilde
P$.

The inclusion $G\hookrightarrow\tilde G$ induces a smooth injection
$G/P\to\tilde G/\tilde P$. Since both spaces have the same dimension,
this must be an open embedding. It is well known that quotients of
semisimple Lie groups by parabolic subgroups are always compact,
whence $G/P\cong\tilde G/\tilde P$. The derivative at the base point
$eP$ of this map is a linear isomorphism $\fg/\frak p\to\tg/\tp$,
which by construction is equivariant over the inclusion
$P\hookrightarrow\tilde P$. Hence the conformal class of inner
products on $\tg/\tp$ from above pulls back to a $P$--invariant
conformal class of inner products on $\frak g/\frak p$. For any Cartan
geometry $(p:\Cal G\to M,\om)$ of type $(G,P)$ we get $TM\cong\Cal
G\x_P\frak g/\frak p$ via $\om$, and hence an induced conformal
structure on $M$.

Nurowski's original construction in \cite{Nurowski} is obtained by
using a degenerate inner product on $\frak g$ to induce an inner
product from the conformal class on $\frak g/\frak p$. Basically, this
amounts to applying the given inner product on $\Bbb R^7$ to the first
columns of matrices. Via the Cartan connection, this is carried over
to a degenerate metric on $\Cal G$, which is shown to induce a well
defined conformal class on $M$. 

\subsection{Weyl structures}\label{3.4} 
We next explain how to obtain the representative metric in Nurowski's
conformal class as described in section \ref{2}. The basic tool is
provided by Weyl structures as introduced in \cite{Weyl}, see also
\cite{CDS} for an alternative approach. For a five manifold $M$ and a
generic rank two distribution $\Cal H\subset TM$, let $p_0:\Cal G_0\to
M$ be the frame bundle of $\Cal H$ and let $(p:\Cal G\to M,\om)$ be
the canonical Cartan geometry. A (local) Weyl structure then is a
$G_0$--equivariant (local) smooth section $\si$ of the natural
projection $\pi:\Cal G\to\Cal G/P_+=\Cal G_0$. There always exist
global Weyl structures, but local ones suffice for our purposes. 

Given a Weyl structure $\si$, one may pull back the Cartan connection
$\om$ to obtain $\si^*\om\in\Om^1(\Cal G_0,\frak g)$. By construction,
for each $i=-3,\dots,3$ the component $\si^*\om_i\in\Om^1(\Cal
G_0,\frak g_i)$ is $G_0$--equivariant. It is better to decompose the
pullback as $\si^*\om=\si^*\om_-+\si^*\om_0+\si^*\om_+$ according to
the decomposition $\fg=\fg_-\oplus\fg_0\oplus\fg_+$. Then the
equivariant form $\si^*\om_-$ descends to an element of
$\Om^1(M,\gr(TM))$. Its value in each point $x\in M$ induces a linear
isomorphism $T_xM\to \gr(T_xM)$ which splits the filtration, i.e.~the
restriction of the $\gr_i(TM)$--component to $T^iM$ coincides with the
canonical projection. Second, the component $\si^*\om_0$ defines a
principal connection on $\Cal G_0$, called the \textit{Weyl
  connection} associated to $\si$. This induces a linear connection on
any vector bundle associated to $\Cal G_0$. Finally, $\si^*\om_+$
descends to a one--form $\Rho\in\Om^1(M,\gr(T^*M))$, called the
\textit{Rho tensor} associated to $\si$. 

Now suppose that we have two tangent vectors $\xi,\eta\in T_xM$. To
compute the values of the metrics in the conformal class on these two
vectors, we have to choose $u\in\Cal G$ with $p(u)=x$ and lifts
$\tilde\xi,\tilde\eta\in T_u\Cal G$. Then we evaluate the elements in
the preferred class of inner products on $\frak g/\frak p$ on
$\om(\tilde\xi)+\frak p$ and $\om(\tilde\eta)+\frak p$. Now we may
linearly identify $\frak g/\frak p$ with $\fg_-$. Denoting elements of
$\fg_-$ as triples $(X,r,Y)$ as suggested by the presentation of
matrices in \ref{3.1}, the conformal class of inner products consists
of all multiples of
\begin{equation}
  \label{confclass}
  ((X,r,Y),(X',r',Y'))\mapsto X^tY'-rr'+Y^tX'.
\end{equation}
Fix a (local) Weyl structure $\si$, choose a point $u_0\in\Cal G_0$
and lifts $\hat\xi,\hat\eta\in T_{u_0}\Cal G_0$ of the two tangent
vectors. Then put $u:=\si(u_0)$, $\tilde\xi=T_{u_0}\si\cdot\hat\xi$
and likewise for $\tilde\eta$. Then the $\fg_-$--components of
$\om(\tilde\xi)$ represent the components of the image of $\xi$ in
$\gr(TM)$ under the isomorphism $TM\cong\gr(TM)$ determined by
$\si$. 

To interpret the individual terms in the right hand side of
\eqref{confclass}, take an element $s\in\fg_2$. Then from the matrix
presentation and the definition of the duality $B$ in \ref{3.1} one
immediately computes that
\begin{gather}
  \label{b1}
  B([s,(X,0,0)],[s,(0,0,Y')])=2s^2X^tY'\\
\label{b2}
B(s,(0,r,0))B(s,(0,r',0))=\frac{s^2}{4}rr'.
\end{gather} 
Passing to associated
bundles, brackets and $B$ correspond to geometric operations. Using
these, one can translate \eqref{confclass} into a geometrically
meaningful formula for a metric. 

Hence it remains to compute the isomorphism $TM\to\gr(TM)$ induced by
some Weyl structure. To pin down one Weyl structure we use an analog
of the scales used in \cite{Weyl}. The associated bundle $\Cal
G_0\x_{G_0}\frak g_2$ is the line bundle $(\gr_{-2}(TM))^*$. For any
Weyl structure $\si$, the Weyl connection $\si^*\om_0$ induces a
linear connection on $(\gr_{-2}(TM))^*$. Since the grading element
acts non--trivially on $\fg_2$, the proof of Theorem 3.8 of
\cite{Weyl} shows that mapping Weyl structures to induced linear
connections on $(\gr_{-2}(TM))^*$ is bijective. In particular, given a
local nowhere vanishing section $\al$ of $(\gr_{-2}(TM))^*$, there is
a unique local Weyl structure such that $\al$ is covariantly constant
for the induced connection. In the language of \ref{2.2} this means
that any generalized contact form for $\Cal H$ determines a Weyl
structure. The main point about the method is that the isomorphism
$TM\to\gr(TM)$ can be computed without knowing the canonical Cartan
connection. In fact, one only has to go through the first steps in the
prolongation/normalization procedure. 

We next describe how to encode the individual parts of a Weyl
structure. Since $\Cal G_0$ is the full frame bundle of $\Cal
H=T^{-1}M$ a principal connection on $\Cal G_0$ is equivalent to a
linear connection $\nabla$ on $T^{-1}M$. Concerning the isomorphism
$TM\to\gr(TM)$, the component in $\gr_{-3}(TM)$ is just given by the
canonical projection $q_{-3}$, so this contains no information.
Suppose that we have given a (local) generalized contact form
$\al\in\Ga((\gr_{-2}(TM))^*)$. Since this is nowhere vanishing, there
is a unique section $\ph\in \Ga(\gr_{-2}(TM))$ such that $\al(\ph)=1$.
Viewing $\al$ as a section of $L(T^{-2}M,\Bbb R)$, the canonical
projection $T^{-2}M\to\gr_{-2}(TM)$ is then given by
$\xi\mapsto\al(\xi)\ph$. Hence we can describe the component in
$\gr_{-2}(TM)$ of the isomorphism $TM\to\gr(TM)$ equivalently by an
extension of $\al$ to a one--form on $M$, which we will again denote
by the same symbol. Finally, the component in $\gr_{-1}(TM)$ of the
isomorphism can be viewed as a projection $\pi_{-1}$ from $TM$ onto
the subbundle $T^{-1}M$. Restricting this projection to $T^{-2}M$, the
kernel is a line subbundle and $q_{-2}$ identifies this line subbundle
with $\gr_{-2}(TM)$. In particular, there is a unique section
$r\in\Ga(T^{-2}M)$ such that $\al(r)=1$ and $\pi_{-1}(r)=0$. 

\subsection{The Weyl structure associated to a generalized contact
  form}\label{3.5} 
Let $\al\in\Ga((\gr_{-2}(TM))^*)$ be a (local) generalized contact
form. As we have seen above, a choice of Weyl structure gives us a
linear connection $\nabla$ on $T^{-1}M$, an extension of $\al$ to a
one--form on $M$, a section $r\in\Ga(T^{-2}M)$, and a projection
$\pi_{-1}$ from $TM$ onto the subbundle $T^{-1}M$. We want to prove
that for the unique Weyl structure such that $\al$ is parallel for the
induced linear connection, these specialize to the objects obtained in
\ref{2.2}--\ref{2.4} (where we used only a part of the connection). We
denote all linear connections induced by our Weyl connection by
$\nabla$. 

The key for verifying this comes from the fact that the Cartan
connection $\om\in\Om^1(\Cal G,\frak g)$ is normal. We have noted in
\ref{3.2} that this implies restrictions on the homogeneity of its
curvature $\Cal K$. For any Weyl structure $\si$, this implies that the
form
\begin{equation}
  \label{Weyl-curv}
  W(\xi,\eta):=d\si^*\om(\xi,\eta)+[\si^*\om(\xi),\si^*\om(\eta)] 
\end{equation}
maps tangent vectors $\xi$ such that $Tp_0\cdot\xi\in T^iM$ and
$Tp_0\cdot\eta\in T^jM$ to $\fg^{i+j+3}$. Now we can split the right
hand side of \eqref{Weyl-curv} into components, which admit a direct
interpretation in terms of the Rho--tensor and a curvature/torsion
quantity $K$ associated to the components of $\si^*\om$, see section 4
of \cite{Weyl}.

For the first step, we will only need components of $K$ with values in
$\gr(TM)$, for which there is an explicit formula in Proposition 4.2
of \cite{Weyl}. For $\ze\in\frak X(M)$ let us denote components in
$\gr(TM)$ under the isomorphism provided by a Weyl structure by
$\ze_i$ for $i=-3,-2,-1$. Then for $\ell<0$, the formula for the
$\fg_\ell$--component of $K(\ze,\ze')$ reads as
\begin{equation}
  \label{K-}
  K_{\ell}(\ze,\ze')=\nabla_\ze\ze'_\ell-\nabla_{\ze'}\ze_\ell
-[\ze,\ze']_\ell+\sum_{i,j<0,i+j=\ell}\{\ze_i,\ze'_j\}. 
\end{equation}

The analysis is best done homogeneity by homogeneity. From Proposition
4.3 of \cite{Weyl} we see that the homogeneous component of degree one
of $K$ coincides with the one of $W$ and hence has to vanish.

\begin{claim}
  Vanishing of the homogeneous component of degree one of $K$ implies
  that $r$ is the Reeb field associated to $\al$ as in 
  \ref{2.2}, the extension of $\al$ to a one--form coincides with the
  one from Proposition \ref{prop23}, and $\nabla$ restricts to the
  partial connection associated to $r$ as in formula
  \eqref{nabla-1-def} in \ref{2.2}.
\end{claim}
\begin{pf}
  The curvature $K$ automatically has positive homogeneity. Hence
  vanishing of the homogeneous component of degree one implies
  vanishing of $K_i$ on $T^{-1}M\x T^iM$ for $i=-3,-2,-1$, and for
  $\xi\in\Ga(T^{-1}M)$ and $\ze\in\Ga(T^iM)$ the value of
  $K_i(\xi,\ze)$ depends only on $q_i(\ze)$. In particular, for $i=2$
  it suffices to compute $K_{-2}(\xi,r)$. Since
  $\al\in\Ga(\gr_{-2}(T^*M))$ is parallel for the induced connection,
  then so is the dual section $\ph=r_{-2}$.  Using that $r_{-1}=0$,
  \eqref{K-} simplifies to give
\begin{equation}
  \label{e2}
  K_{-2}(\xi,r)=\al([\xi,r])\ph=[\xi,r]_{-2}.  
\end{equation}
Next, for $\xi\in\Ga(T^{-1}M)$ and $\ze\in\frak X(M)$, we
obtain from \eqref{K-}
\begin{equation}
  \label{e3}
  K_{-3}(\xi,\ze)=\nabla_\xi \ze_{-3}-[\xi,\ze]_{-3}+\{\xi,\ze_{-2}\}. 
\end{equation}
Now put $\ze:=[r,\eta]$ for $\eta\in\Ga(T^{-1}M)$. Then
$\ze_{-3}=q_{-3}(\ze)=\{\ph,\eta\}$, and since $\nabla$ is compatible
with $\{\ ,\ \}$ and $\nabla\ph=0$ we see that 
\begin{equation}
  \label{e4}
 K_{-3}(\xi,[r,\eta])=\{\ph,\nabla_\xi\eta\}-q_{-3}([\xi,[r,\eta]])-
\{\xi,K_{-2}(\eta,r)\}. 
\end{equation}
Vanishing of \eqref{e3} and \eqref{e2} thus implies that $\nabla$ is
the connection determined by $r$. But then the fact that $\ph$ is
parallel implies that $r$ is the Reeb field associated to $\al$ as in
\ref{2.2}.  Given this, vanishing of \eqref{e2} says that
we get the right extension of $\al$ to a one--form.\qed 
\end{pf}

\begin{rem}
  It can be actually shown that the opposite implication holds as
  well. If we use $r$, the extension of $\al$ and $\nabla$ as the data
  associated to the Weyl form and $\al$ is parallel for the induced
  connection, then the facts that $r$ is the Reeb field, we have the
  right extension of $\al$, and $\nabla$ restricts to the partial
  connection determined by $r$ imply that the homogeneous component of
  degree one of the curvature $K$ vanishes.
\end{rem}

It remains to show that our Weyl structure produces the right
projection $\pi_{-1}$. For this we have to analyze the homogeneous
components of degree two of $W$ and $K$. According to Proposition 4.3
of \cite{Weyl}, the difference between these two components is
determined the homogeneous component of degree 2 of the Rho--tensor.
Since we will not need any other part of the Rho--tensor, we simply
denote this component by $\Rho$. It can be either interpreted as a
partially (on $T^{-1}M$) defined one--form with values in
$(T^{-1}M)^*$ or as a bilinear form on $T^{-1}M$. Further, we will
also need components of $W$ and $K$ in degree zero. These are sections
of the bundle $\Cal G_0\x_{G_0}\fg_0$, so in particular such a section
induces an endomorphism of $\gr_i(TM)$ for $i=-3,-2,-1$. This action
is induced by the components $\fg_0\x\fg_i\to\fg_i$ of the Lie bracket
on $\fg$.

We will also need some of the other tensorial maps induced by the Lie
brackets on $\fg$, and we will denote them all by $\{\ ,\ \}$. In
particular, these define bilinear bundle maps
$\gr_{-i}(TM)\x\gr_i(T^*M)\to \Cal G_0\x_{G_0}\fg_0$ for all $i$, as
well as 
$$\gr_{-j}(TM)\x\gr_i(T^*M)\to\gr_{i-j}(TM)
$$ 
for $i<j$. 

For $K$, the component $K_0$ is the curvature of the Weyl connection,
see Proposition 4.2 of \cite{Weyl}. Hence the induced
endomorphism on $\gr_i(TM)$ is simply the curvature $R$ of the
corresponding linear connection. Since we need the component of degree
2, we are interested in $K_0$ and $W_0$ as two--forms acting on
$T^{-1}M\x T^{-1}M$, and there the difference between $W$ and $K$ is
given by 
$$
(\xi,\eta)\mapsto \{\Rho(\xi),\eta\}-\{\Rho(\eta),\xi\}.
$$

Now we first observe that the bundle $\gr_{-2}(TM)$ admits the nonzero
parallel section $\phi$, so $R$ has to act trivially on
$\gr_{-2}(TM)$. Hence vanishing of the restriction of $W_0$ to
$T^{-1}M\x T^{-1}M$ implies that also
$\{\Rho(\xi),\eta\}-\{\Rho(\eta),\xi\}$ acts trivially on
$\gr_{-2}(TM)$. But one immediately verifies that for $Z\in\fg_1$
and $X\in\fg_{-1}$ the action of $[Z,X]\in\fg_0$ on $\fg_{-2}$ is by
multiplication by a nonzero multiple of $B(Z,X)=ZX$. Hence we conclude
that, viewed as a bilinear form on $T^{-1}M$, $\Rho$ is symmetric.

Further, one verifies directly that for $Z\in\fg_1$ and
$X_1,X_2\in\fg_{-1}$, one has 
\begin{align*}
  [[Z,X_1],X_2]&=B(Z,X_1)X_2-3B(Z,X_2)X_1\\
  [Z,[X_1,X_2]]&=4(B(Z,X_1)X_2-B(Z,X_2)X_1).
\end{align*}
Using these two identities and the symmetry of $\Rho$ one immediately
verifies that
$$
\{\{\Rho(\xi),\eta\}-\{\Rho(\eta),\xi\},\xi'\}=
-\tfrac{3}{4}\{\Rho(\xi'),\{\xi,\eta\}\}
$$
for all $\xi,\eta,\xi'\in T^{-1}M$. 

Now let us assume that $\{\xi,\eta\}=\ph$. Then by definitions of the
curvature $R$ and of $\Ph$ in \eqref{Phidef} in \ref{2.4} we get
$R(\xi,\eta)(\xi')=\Ph(\xi')-\nabla_r\xi'$. Hence we conclude that
vanishing of $W_0(\xi,\eta)$ implies (renaming $\xi'$ to $\xi$) that
\begin{equation}
  \label{e5}
 \Ph(\xi)-\nabla_r\xi+\tfrac{3}{4}\{\ph,\Rho(\xi)\}=0
\end{equation}
for all $\xi\in\Ga(T^{-1}M)$. 

For the remaining components, we can use formula \eqref{K-} from
\ref{3.5} to compute $K$, and the correction to $W$
is given by those $\Rho$--terms which involve entries from
$T^{-1}M$. Vanishing of $W_{-1}(r,\xi)$ for $\xi\in\Ga(T^{-1}M)$ implies 
\begin{equation}
  \label{e6}
  \nabla_r\xi-\pi_{-1}([r,\xi])+\{\ph,\Rho(\xi)\}=0.
\end{equation}
Finally, vanishing of $W_{-3}(r,\ze)$ for $\ze\in\frak X(M)$ gives
$$
\nabla_rq_{-3}(\ze)-q_{-3}([r,\ze])+\{\ph,\pi_{-1}(\ze)\}=0.
$$
Inserting $\ze=[r,\xi]$ and using equation \eqref{Psidef} from
\ref{2.4} we see that we can pull off $\{\ph,\ \}$ to conclude that 
\begin{equation}
  \label{e7}
  \nabla_r\xi-2\Ps(\xi)+\pi_{-1}([r,\xi])=0.
\end{equation}
Using \eqref{e6} to compute $\{\ph,\Rho(\xi)\}$ and \eqref{e7} to
compute $\nabla_r\xi$ and inserting both into \eqref{e5}, we obtain 
$$
\pi_{-1}([r,\xi])-\tfrac{2}{5}\Ph(\xi)+\tfrac{7}{5}\Ps(\xi)=0,
$$
which exactly means that we get the right projection $\pi_{-1}$. 

\begin{rem} 
  The formula for $[Z,[X_1,X_2]]$ from above shows that the bracket
  $\{\ ,\ \}:\gr_1(T^*M)\x\gr_{-2}(TM)\to\gr_{-1}TM$ is explicitly
  given by $\{\ps,\{\xi,\eta\}\}=4(\ps(\xi)\eta-\ps(\eta)\xi)$. The
  other components of $\{\ ,\ \}$ can be computed similarly. Using
  these formulae, one easily verifies that the transformation laws in
  Lemmas \ref{lem25} and \ref{lem26} are the specializations of Proposition
  3.4. of \cite{Weyl}, which gives a general formula for the change of
  the isomorphism $TM\to\gr(TM)$ caused by a change of Weyl structure.
\end{rem}

\subsection{Computing the metric}\label{3.7}
With the description of the isomorphism $TM\to\gr(TM)$ at hand, we can
now verify the formula for the metric. We only have to interpret the
expressions \eqref{b1} and \eqref{b2} from \ref{3.4} in geometric
terms. In these formulae, $s\in\fg_2$ corresponds to the generalized
contact form $\al$. Let us further suppose that $\ze,\ze'\in\frak
X(M)$, are vector fields. Then in \eqref{b2}, the element
$r\in\fg_{-2}$ corresponds to $\ze_{-2}=\al(\ze)\ph$ and likewise for
$r'$. Thus the geometric interpretation of \eqref{b2} simply is
$\al(\ze)\al(\ze')$. 

In \eqref{b1}, the element $X\in\fg_{-1}$ corresponds to
$\pi_{-1}(\ze)$ and $Y'$ corresponds to
$q_{-3}(\ze')=\{\ph,\ze'_1\}$. Hence what we actually have to do is
interpreting (again in the notation of \ref{3.4})
$$
B([s,(X,0,0)],[s,[r_0,(X',0,0)]]),
$$ 
where $r_0\in\fg_{-2}$ is characterized by $B(r_0,s)=1$. Now one
easily computes that in the Lie algebra $\fg$, one has
$[s,[r_0,X']]=3X'$ and using this, one verifies that 
$B([s,X],[s,[r_0,X']])r_0=6[X,X']$. Using that 
$$
\{\xi,\eta\}=\al([\xi,\eta])\ph=-d\al(\xi,\eta)\ph, 
$$
we see that the geometric interpretation of \eqref{b1} is
$-6d\al(\pi_{-1}(\ze),\ze'_1)$. But then formula \eqref{confclass}
from \ref{3.4} shows that Nurowski's conformal class contains the
metric
$$
(\ze,\ze')\mapsto
-3d\al(\pi_{-1}(\ze),\ze'_1)-3d\al(\pi_{-1}(\ze'),\ze_1)-4\al(\ze)\al(\ze'), 
$$ 
which proves
\begin{thm}\label{thm37}
  The metric $g_\al$ defined in formula \eqref{metric} in \ref{2.4} is
  contained in Nurowski's conformal class. 
\end{thm}


\begin{thebibliography}{XX}

\bibitem{BDS} D. Burns, K. Diederich, S. Shnider, Distinguished curves
  in pseudoconvex boundaries, Duke Math.  J. \textbf{44} no. 2 (1977)
  407--431.

\bibitem{CDS} D.M.J. Calderbank, T. Diemer, V. Sou\v cek,
  Ricci--corrected derivatives and invariant differential operators.
  Differential Geom. Appl. \textbf{23}, no. 2 (2005) 149--175. 

\bibitem{Srni05} A. \v Cap, Two constructions with parabolic
  geometries, Rend. Circ. Mat. Palermo Suppl. ser. II, \textbf{79}
  (2006) 11--37. 

\bibitem{Weyl} A. \v Cap, J. Slov\'ak, Weyl structures for parabolic
  geometries, Math. Scand. \textbf{93}, No.~1 (2003) 53-90.

\bibitem{Cartan:five} E. Cartan, Les syst\`emes de Pfaff a cinq variables
et les \'equations aux deriv\'ees partielles du second ordre,
Ann. Ec. Normale \textbf{27} (1910), 109--192.

\bibitem{Nurowski} P. Nurowski, Differential equations and conformal
  structures,  J. Geom. Phys. \textbf{55}, no. 1 (2005) 19--49.

\bibitem{Lee1} J.M. Lee, \textit{The Fefferman metric and
    pseudo-Hermitian invariants}, Trans. Amer. Math. Soc.
  \textbf{296}, no. 1 (1986) 411--429.

\bibitem{Lee2} J.M. Lee, \textit{Pseudo-Einstein structures on CR
    manifolds}, Amer. J. Math. \textbf{110}, no. 1 (1988) 157--178.

\bibitem{Sa1} K. Sagerschnig, Parabolic geometries determined by
  filtrations of the tangent bundle, Rend. Circ. Mat. Palermo Suppl.
  ser. II, \textbf{79} (2006) 175--181.

\bibitem{Sa2} K. Sagerschnig, Split octonions and generic rank two
  distributions in dimension five, Arch. Math. (Brno) \textbf{42}
  Suppl. 329--339.

\bibitem{Springer} T.A. Springer, F.D. Feldenkamp, Octonions, Jordan
  algebras and exceptional groups. Springer, Berlin, 2000. 

\bibitem{Tanaka} N. Tanaka, On the equivalence problem associated
  with simple graded Lie algebras. Hokkaido Math. J., \textbf{8}
  (1979), 23--84.

\end{thebibliography}
\end{document}